\documentstyle[newlfont,amstex,amssymb]{article}
\newtheorem{theorem}{\mediumseries\sc Theorem}
\newtheorem{lemma}[theorem]{\mediumseries\sc Lemma}

\newtheorem{definition}[theorem]{\mediumseries\sc Definition}

\title{Donder's Version of Revised Countable Support}
\author{Ulrich Fuchs\thanks{{\bf\sf Freie Universit\"at Berlin,}
e-mail: {\bf\tt FuX\protect@\noexpand@Math.FU-Berlin.De}}}

\begin{document}
\maketitle
\begin{abstract}
Shelah introduced the revised countable support (RCS) iteration to iterate
semiproperness~\cite{shelah}.
This was an endpoint in the search for an iteration of a weak condition,
still implying that $\boldsymbol\aleph\sb1$ is preserved.
It was one of the key tools in the proof of the relative consistency of
{\sf Martin's Maximum}.
Dieter Donder found a better manageable approach to this iteration, which is
presented here.
More iterations of semiproperness are formulated in~\cite{schlindwein}
and~\cite{miyamoto}.
\end{abstract}
\section*{Terminology}
An {\bf iteration} is a well-ordered commuting family of complete embeddings
of complete Boolean algebras.
From now on let
$\cal B=\left\langle\Bbb B\sb\gamma|\gamma\in\lambda\right\rangle$
be an iteration, which hence means
$\forall\alpha<\beta<\lambda\quad\Bbb B\sb\alpha\lessdot\Bbb B\sb\beta$.
We will formulate in this section statements and definitions of
notions for this
special case only, but they all have formulations and generalizations
for the general case.

If $\alpha<\lambda$ we have $\Vdash\sb\alpha(\text{$\cal B/G$ is iteration})$.
(A filter $G$ in $\Bbb B\sb\alpha$ defines
quotients $\Bbb B\sb\beta/G$ $(\beta>\alpha)$ and embeddings
$\Bbb B\sb\beta/G\longrightarrow\Bbb B\sb\gamma/G$
$(\gamma>\beta>\alpha)$ even if the algebras in $\cal B$ are not complete.
The system is called $\cal B/G$).

Of course the family of the canonical projections
$h\sb\alpha:\bigcup\cal B\longrightarrow\Bbb B\sb\alpha$,
$b\mapsto\bigwedge\{c\in\Bbb B\sb\alpha| b\leqslant c\}$
does also commute.
(i.e.~$\alpha<\beta\Longrightarrow h\sb\alpha\circ h\sb\beta=h\sb\alpha$).

A {\bf thread} in $\cal B$ is a $f\in\prod\cal B$ with
$\forall\alpha\leqslant\beta<\lambda\quad f(\alpha)=h\sb\alpha(f(\beta))$.

If $f$ is a thread, we can define an iteration $\cal B{\restriction} f$ as
$\left\langle\Bbb B\sb\alpha{\restriction}
f(\alpha)|\alpha<\lambda\right\rangle$
together with the complete embeddings
$\Bbb B\sb\alpha{\restriction} f(\alpha)\longrightarrow\Bbb B\sb\beta
{\restriction} f(\beta)$, $b\mapsto b\wedge f(\beta)$.

$\cal T(\cal B)$ is the set of all threads, which is canonically
componentwise partially ordered.

Let $c:\bigcup\cal B\longrightarrow\cal T(\cal B)$ be defined by
$b\mapsto\left\langle h\sb\alpha(b)|\alpha<\lambda\right\rangle$.
The range of $c$ is the set $\cal C(\cal B)$
of all eventually constant threads.

Let $\cal C(\cal B)\subset F\subset\cal T(\cal B)$.
If
\begin{equation}\label{threads}
\forall\alpha<\lambda\quad\forall f\in F\quad
\forall b\in\Bbb B\sb\alpha\quad\Bigl(b\leqslant f(\alpha)\Longrightarrow
\quad\text{$f$ and $c(b)$ are compatible in $F$}
\Bigr)
\end{equation}
then $F$ is separative and the mappings
$\Bbb B\sb\alpha\longrightarrow F$, $b\mapsto c(b)$ are complete embeddings.
Hence there is a complete Boolean algebra $\Bbb B(F)$ and a dense embedding
$d:F\longrightarrow\Bbb B(F)$ such that $d\circ c\subset\text{\sf id}$
i.e.~$\stackrel\frown{\cal B\Bbb B}(F)$ is an iteration.

$\cal C(\cal B)$ and $\cal T(\cal B)$ satisfy (\ref{threads}).~Let
$\operatorname{\text{\normalshape\sf Dir}}(\cal B)=\Bbb B(\cal C(\cal B))$
and
$\operatorname{\text{\normalshape\sf Inv}}(\cal B)=\Bbb B(\cal T(\cal B))$.

{\sc Factor Property:} Let
$\Bbb E=\operatorname{\text{\normalshape\sf Inv}}(\cal B)$
and
$\Bbb D=\operatorname{\text{\normalshape\sf Dir}}(\cal B)$.
For all $\alpha<\lambda$
\[
\Vdash\sb\alpha
\bigl(
\text{$\Bbb E/G\cong\operatorname{\text{\normalshape\sf Inv}}({\cal B/G})$
\quad and\quad
$\Bbb D/G\cong\operatorname{\text{\normalshape\sf Dir}}({\cal B/G})$}
\bigr).
\]
{\sc Fact:} (proof in \cite{jech}, Lemma~36.5, page~460).
If $\lambda>\omega$ is regular and
$\left\{\alpha\mid\Bbb B_\alpha=
\operatorname{\text{\normalshape\sf Dir}}(\cal B{\restriction}\alpha)\right\}$
is stationary in $\lambda$ and each $\Bbb B_\alpha$ satisfies the
$\lambda$-antichain-condition, then so does
$\operatorname{\text{\normalshape\sf Dir}}(\cal B)$.
\section*{Revised countable support iterating semiproperness}
Let us assume familiarity with the following facts about the iteration
of semiproperness, more or less proved by simultaneously playing
semiproper games in different generic extensions:
\begin{theorem}\label{theorem}
\mbox{}
\begin{enumerate}
\item
Let $\Bbb P$ be semiproper and
$\Vdash\sb{\Bbb P}(\text{$\dot{\Bbb Q}$ is semiproper})$,
then $\Bbb P*\Bbb Q$ is semiproper.
\item\label{seminv}
Let $\cal B=\langle\Bbb B\sb n| n\in\omega\rangle$
be some iteration with semiproper $\Bbb B\sb0$ and
$\forall n\quad\Vdash\sb n(\text{$\Bbb B\sb{n+1}/G$ is semiproper})$.
Then $\operatorname{\text{\normalshape\sf Inv}}(\cal B)$ is semiproper.
\item\label{semidir}
Let $\lambda>\omega$ be regular and
$\cal B=\langle\Bbb B\sb\alpha|\alpha\in\lambda\rangle$
be an iteration with: for all $\alpha<\lambda$
$\Bbb B\sb\alpha$ is semiproper and
\begin{gather*}
\forall\beta\in(\alpha,\lambda)
\quad\Vdash\sb\alpha(\text{$\Bbb B\sb\beta/G$ is semiproper})
\\
\operatorname{\text{\normalshape\sf cof}}(\alpha)=\omega
\quad\Longrightarrow\quad\Bbb B\sb\alpha=
\operatorname{\text{\normalshape\sf Inv}}(\cal B{\restriction}\alpha)
\end{gather*}
If $\lambda=\omega\sb1$ or $\operatorname{\text{\normalshape\sf Dir}}(\cal B)$
satisfies the $\lambda$-antichain-condition, then
$\operatorname{\text{\normalshape\sf Dir}}(\cal B)$ will be semiproper.
\end{enumerate}
\rule{4cm}{0.3pt}
\end{theorem}
\begin{definition}
We call a thread $f\in\cal T(\cal B)$ {\bf short} iff
\[
\exists\alpha<\lambda\quad f(\alpha)\Vdash\sb\alpha
(\operatorname{\text{\normalshape\sf cof}}(\lambda)=\omega).
\]
$\cal S(\cal B)$, the set of all short or eventually constant threads,
satisfies (\ref{threads}). Define
$\operatorname{\text{\normalshape\sf Rlim}}(\cal B)=\Bbb B(\cal S(\cal B))$.
\end{definition}
We remark that {\sf Rlim} satsfies the factor property.
Now our RCS-iteration, that means taking {\sf Rlim} at limit stages,
of semiproperness will work:
\begin{theorem}
Let
$\cal B=\left\langle{{\Bbb B}\sb\alpha}|\alpha\leqslant\lambda\right\rangle$
be a RCS-iteration satisfying%
\begin{align*}
&\text{$\Bbb B\sb0$ is semiproper}\\
&\forall\alpha<\lambda\quad
\begin{aligned}[t]
&\Vdash\sb{\alpha+1}(|\Bbb B\sb\alpha|\leqslant\boldsymbol\aleph\sb1)\\
&\Vdash\sb\alpha(\text{$\Bbb B\sb{\alpha+1}/G$ is semiproper}).
\end{aligned}
\end{align*}
Then $\Bbb B\sb\lambda$ is semiproper.%
\end{theorem}
Prove by induction on $\beta\leqslant\lambda$ that
$\forall\alpha<\beta\quad
\Vdash\sb\alpha(\text{$\Bbb B\sb\beta/G$ is semiproper})$.
In the limit case use the factor property of {\sf Rlim} and the
following lemma:
\begin{lemma}\label{lemma}
Let $\cal B=\langle{\Bbb B\sb\alpha}|\alpha<\lambda\rangle$
be RCS-iteration, $\lambda$ limit ordinal such that:
\begin{align*}
\forall\alpha<\lambda\quad
&\text{$\Bbb B\sb\alpha$ is semiproper}
\\
&\Vdash\sb{\alpha+1}(\left|\Bbb
B\sb\alpha\right|\leqslant\boldsymbol\aleph\sb1)
\\
&\forall\beta\in(\alpha,\lambda)\quad
\Vdash\sb\alpha(\text{$\Bbb B\sb\beta/G$ is semiproper}).
\end{align*}
Then $\operatorname{\text{\normalshape\sf Rlim}}(\cal B)$ is semiproper.
\end{lemma}
{\sc Proof of lemma~\ref{lemma}:}
Let $\Bbb B=\operatorname{\text{\normalshape\sf Rlim}}(\cal B)$ and
$d:\cal S(\cal B)\longrightarrow\Bbb B$ be dense.
Let us show that
$\left\{b\bigm|\text{$\Bbb B{\restriction} b$ is semiproper}\right\}$
is dense in $\Bbb B$.
So let $f\in\cal S(\cal B)$ and we have to find a $g\in\cal S(\cal B)$
below $f$ such that $\Bbb B{\restriction} d(g)$ is semiproper.
\begin{itemize}
\item
first case:
$\exists\alpha<\lambda\quad f(\alpha)\not\Vdash\sb\alpha
(\operatorname{\text{\normalshape\sf cof}}(\lambda)>\omega)$.

Hence there is a thread $g$ below $f$ and an $\alpha<\lambda$ such that
$g(\alpha)\Vdash\sb\alpha
(\operatorname{\text{\normalshape\sf cof}}(\lambda)=\omega)$.
We have
$\Bbb B{\restriction} d(g)\cong
\operatorname{\text{\normalshape\sf Inv}}(\cal B{\restriction} g)$
(every thread below $g$ is short).
In a generic extension over
$\Bbb B\sb\alpha{\restriction} g(\alpha)$
we apply (up to isomorphism)~\ref{theorem}.\ref{seminv} and get
\[
\Vdash\sb{\Bbb B\sb\alpha{\restriction} g(\alpha)}
(\text{$\operatorname{\text{\normalshape\sf Inv}}(\cal B{\restriction} g/G)$
is semiproper}),
\]
the factor property for the inverse limit gives
$\Bbb B{\restriction} d(g)$ is semiproper.
\item
second case:
$\forall\alpha<\lambda\;f(\alpha)\Vdash\sb\alpha
(\operatorname{\text{\normalshape\sf cof}}(\lambda)>\omega)$.

Hence below $f$ there is no short thread.
We distinguish two more cases:
\begin{itemize}
\item
$\exists\alpha<\lambda\quad f(\alpha)\not\Vdash\sb\alpha
(\operatorname{\text{\normalshape\sf
cof}}(\lambda)>\boldsymbol\aleph\sb1)$.

Hence there is an $\alpha<\lambda$ and an eventually constant thread $g$
below $f$ such that
$g(\alpha)\Vdash\sb\alpha
(\operatorname{\text{\normalshape\sf cof}}(\lambda)=\omega\sb1)$.
Now proceed as in the first case, using~\ref{theorem}.\ref{semidir} instead
of~\ref{theorem}.\ref{seminv} and the factor property for the direct limit
to get
$\Bbb B{\restriction} d(g)\cong
\operatorname{\text{\normalshape\sf Dir}}(\cal B{\restriction} g)$ is semiproper.
\item
it remains:
$\forall\alpha<\lambda\;f(\alpha)\Vdash\sb\alpha
(\operatorname{\text{\normalshape\sf cof}}(\lambda)>\boldsymbol\aleph\sb1)$.

So $\forall\alpha<\lambda\quad|B\sb\alpha|<\lambda$,
and $\lambda$ is regular. (Otherwise the cofinality would have been collapsed).
Let
$S=\left\{\alpha<\lambda|\operatorname{\text{\normalshape\sf cof}}(\alpha)=
\boldsymbol\aleph\sb1\right\}$.
We have
$\forall\alpha\in S\quad\Vdash\sb\alpha
(\operatorname{\text{\normalshape\sf cof}}(\alpha)=\boldsymbol\aleph\sb1)$
and therefore
$\Bbb B\sb\alpha{\restriction} f(\alpha)\cong
\operatorname{\text{\normalshape\sf Dir}}
((\cal B{\restriction} f){\restriction} \alpha)$.
Since $S$ is stationary in $\lambda$,
$\operatorname{\text{\normalshape\sf Dir}}(\cal B{\restriction} f)$
satisfies the $\lambda$-antichain-condition.
(apply the {\sc Fact} of the previous section).

Now we can apply~\ref{theorem}.\ref{semidir} and get the semiproperness
of $\Bbb B{\restriction} d(f)\cong
\operatorname{\text{\normalshape\sf Dir}}(\cal B{\restriction} f)$.
\end{itemize}
\end{itemize}
\rule{4cm}{0.3pt}

\end{document}